\documentclass{article}

\usepackage{arxiv}

\usepackage{natbib}

\usepackage{algorithm}
\usepackage{algorithmic}

\usepackage{hyperref}
\usepackage{url}

\title{Finding Optimal Policy for Queueing Models:\\ New Parameterization}

\author{Trang H. Tran$^{1}$, Lam M. Nguyen$^{2}$, \textbf{Katya Scheinberg}$^{1}$ \\
$^{1}$ School of Operations Research and Information Engineering, Cornell University, Ithaca, NY, USA\\
$^{2}$ IBM Research, Thomas J. Watson Research Center, Yorktown Heights, NY, USA\\
\texttt{htt27@cornell.edu},
\texttt{LamNguyen.MLTD@ibm.com}, \texttt{katyas@cornell.edu}
}

\usepackage{pifont}

\usepackage{epsfig}
\usepackage{amssymb}
\usepackage{amsmath}
\usepackage{amsthm}
\usepackage{amsfonts}
\usepackage{bbding}
\usepackage{array,listings,xcolor}
\usepackage{paralist}
\usepackage{xargs}                      
\usepackage{caption}

\usepackage{algorithm}
\usepackage{algorithmic}

\usepackage{hyperref}
\usepackage{url}

\newtheorem{thm}{Theorem}

\newcolumntype{C}[1]{>{\centering\let\newline\\\arraybackslash\hspace{0pt}}m{#1}}

\newcommand{\E}{\mathbb{E}} 
\newcommand{\N}{\mathbb{N}}

\usepackage{xcolor}

\newcommand{\R}{\mathbb{R}}

\newcommand{\argmax}{\mathrm{arg}\!\displaystyle\max}

\newcommand{\zero}[1]{{\boldsymbol{0}}}

\DeclareMathOperator{\A}{\mathcal{A}} 
\renewcommand{\S}{\mathcal{S}}

\definecolor{codegreen}{rgb}{0,0.6,0}
\definecolor{codegray}{rgb}{0.5,0.5,0.5}
\definecolor{codepurple}{rgb}{0.58,0,0.82}
\definecolor{backcolour}{rgb}{0.98,0.98,0.98}

\lstdefinestyle{mystyle}{
    language = Python,
    backgroundcolor=\color{backcolour},   
    commentstyle=\color{codegreen},
    keywordstyle=\color{magenta},
    numberstyle=\tiny\color{codegray},
    stringstyle=\color{codepurple},
    basicstyle=\ttfamily\footnotesize,
    breakatwhitespace=false,         
    breaklines=true,                 
    captionpos=b,                    
    keepspaces=true,                 
    showspaces=false,                
    showstringspaces=false,
    showtabs=false,                  
    tabsize=2
}
\hypersetup{
  colorlinks   = true, 
  urlcolor     = blue, 
  linkcolor    = blue, 
  citecolor   = blue 
}

\begin{document}

\maketitle

\begin{abstract}
Queueing systems appear in many important real-life applications including communication networks, transportation and manufacturing systems. Reinforcement learning (RL) framework is a suitable model for the queueing control problem where the underlying dynamics are usually unknown and the agent receives little information from the environment to navigate. In this work, we investigate the optimization aspects of the queueing model as a RL environment and provide insight to learn the optimal policy efficiently. We propose a new parameterization of the policy by using the intrinsic properties of queueing network systems. Experiments show good performance of our methods with various load conditions from light to heavy traffic.
\end{abstract}


\section{Introduction}\label{sec_intro}

The optimal control problem of queueing networks has many important applications in real life, including communications networks, transportation, and manufacturing systems \citep{reentrant, facturing, communication,Nguyen2017_fluid, inpatient_jim}.
The main goal is to efficiently determine the optimal control policy, especially when the network is heavily loaded. In certain  cases, such a policy can be derived analytically as a function of system parameters, but even in those simple cases, since system parameters are rarely known in practice, applying a learning approach is often preferable \citep{ieee_baixie2019}. Queueing networks can be modeled effectively as a Reinforcement Learning (RL) environment, where an agent chooses actions in reaction to the environment and receives feedback to learn a good policy. RL has been a very active area recently, with the success of well-known algorithms such as AlphaGo \citep{Silver2016} and OpenAIFive \citep{OpenAI_dota} and is expected to have tremendous impact on real-life applications, including autonomous driving, dynamic treatment regimes in healthcare and robotics manipulation. Despite its practical success, there is not enough understanding of the algorithmic frameworks used in this rich environment, especially, the optimization process. 

The main goal of this paper is to investigate the key steps of the optimization process in model-free RL using a particular practical and yet tractable application of simple queuing network control. There have been recent works applying Reinforcement Learning framework to the queueing systems control \citep{ieee_baixie2019, dai2020queueing}. However, these papers use a model-based approach where the agent estimates the system parameters in order to create a model of the system. In this paper, we are interested in  investigating performance of the model-free framework on this problem, 
 without  explicitly modeling and utilizing the system dynamics.

We choose to consider a simple parallel queueing system with one server, exponential interarrival and service times. The optimal policy for this system is a well known threshold-type and can be computed based on the service time - holding cost ratio of each queue. A RL method is expected to learn this policy without the explicit knowledge of the service and arrival rates and possibly even without knowledge of the holding costs by simply maximizing a reward function (or minimizing the cost function), which is computed by simulating the queueing system over a set of parameters that define the policy. There are many choices which influence the efficiency of the optimization process and the effects of most of these choices are poorly explored in the literature.  For example, while the classical REINFORCE estimator \citep{Williams1992} has been used widely in many RL applications, it is not clear how it behaves in comparisons with other gradient estimators, and in particular,  on queueing models. In addition, the choice of reward functions and policy parameterizations have a significant effect on the nature of the resulting optimization problem and thus on the optimization process.
Here we investigate this effect by considering some alternative choices for model-free RL applied to queueing systems.

\textbf{Contributions:}
\begin{itemize}
    \item We investigate a simple linear policy parametrization which can be used within a model-free RL framework and show that it is able to approximate the optimal threshold with arbitrary accuracy. Based on the representation and approximation of the (optimal) priority policy, we then propose a new logarithm-scale parameterization which changes the problem scale and allows for smoother simpler optimization landscapes.  

    \item We compared two different methods of constructing gradient approximation: the REINFORCE estimator \cite{Williams1992} and the Finite Difference gradient estimator \citep{fd_pde,aaai_fd2}. 
    We implement the two gradient estimators with an adaptive line-search algorithm which is robust to gradient estimators and whose performance can be easily derived as a function of the cost and the error of the gradient estimators. 
    \item Finally, we show that logarithm-scale parameterization behaves similarly or better than linear-scale parameterization regardless of the choice of the gradient estimator. 
    \end{itemize}

\textbf{Related Works.}

\textbf{\textit{Reinforcement Learning.}}
Reinforcement learning techniques have gained significant  popularity recently for a large variety of applications, while the general  RL problem setting remains difficult without specificity of the environment and the  structure of the objective \citep{bookrl,Sutton2018}. 
The  framework is generally described as a Markov decision process (MDP), with a particular transition probability function and a reward mechanism that returns the feedback after the agent takes an action. The typical algorithm for solving an MDP is via policy iteration or value iteration \citep{valueiteration}, however, this approach is often not suitable in learning applications when the dimension of the problem is large. 
Another approach is to estimate the state-value function or the value function such as Q-learning \citep{Watkins1992} and its variants  (e.g. Deep Q-learning (DQN) \citep{Mnih2013, Mnih2015}, and double Q-learning \citep{Hasselt2016}). However, it has been observed that learning the state-value function is not efficient when the action space is large or even infinite.  

\textbf{\textit{Parameterized Policy and Optimization.}}
More popular recent approach to RL is to learn the policy directly as a parameterized function which maps  the state space to a distribution over the action space. Neural networks are typical candidates for such  parameterized functions because they have shown efficiency in many machine learning applications \citep{Sutton2018,bookrl,dai2020queueing}. However, in this paper we show that a simpler linear parameterization is sufficient to model the class of policies of the particular queueing system under consideration. 

After  parameterization of the control policy is chosen, various optimization methods can be applied to find the optimal value of the parameters with respect to the expected reward function.  Usually such a method involves some version of a gradient ascent using estimates of the gradient of the reward function with respect to the parameters. 
One of the classical approaches is REINFORCE \citep{Williams1992}, which computes an estimator of the policy gradient using likelihood ratio \citep{stoc_opt, glynn}. However, the likelihood ratio estimator is known to have high variance and REINFORCE estimator also suffers from that weakness. There have been other improvements to reduce the variance e.g. adding baseline terms   or discarding some rewards in the GPOMDP estimator \citep{Sutton2018,Zhao2011,Baxter2001, pmlr-v108-pham20a} . 
Advanced policy optimization algorithms include Trust Region Policy Optimization (TRPO) \citep{Schulman2015} and Proximal Policy Optimization \citep{Schulman2017} which go beyond simple gradient accent with the aim to better control the steps and progress of the optimization algorithm. 

The update step of the REINFORCE algorithm has the foundation from stochastic gradient descent (SGD) method \citep{Robbins1951}. 
SGD with its stochastic first-order variants \citep{AdaGrad,Kingma2014,Bottou2018,Nguyen2018_sgdhogwild} and variance reduction methods \citep{SAG, SAGA, SVRG,Nguyen2017sarah} are favorable in machine learning because they are efficient in dealing with large-scale problems. 
In fact, other gradient estimators can be used to update the policy in place of the likelihood ratio method. Finite Difference estimators have been one of the popular choices, which is applied successfully to numerous applications in robotics \citep{pg_robotics}. It has been observed that finite difference (FD) estimators generally have lower variance than the likelihood ratio estimator in typical REINFORCE algorithm \citep{pg_robotics}, which may result in a faster convergence of the optimization algorithm. However, such FD estimators
are more computationally costly when the problem dimension is large. 
To properly compare both type of gradient estimators in application to queuing systems we embed them within an adaptive first order optimization method  designed for noisy stochastic gradient and function estimators \citep{katya_linesearch,billy_miaolan_katya}. 
This method is simpler than methods for policy gradients (e.g. TRPO and PPO \citep{Schulman2015,Schulman2017}) but is more advanced  than those using FD estimators. Applying this method for optimization when the optimal policy is known 
helps us gain knowledge about different gradient estimators and compare their performance.

\textbf{\textit{Queueing systems control problem.}}
We consider queueing networks and stochastic processing networks, where ``job" arrive into the system and join a queue to be processed by a server \citep{reentrant, facturing, communication,Nguyen2016_agent, inpatient_jim}.  The control problem is defined as a process of allocating servers to processing jobs. In the queuing networks where interarrival times and service times are exponentially distributed, the control problem can be modeled as a Markov decision process (MDP) via uniformaization \citep{unif_cont_discrete, mdp_uniformization}. The objective function of interest is often the holding cost, where jobs in queue $i$ incur a holding cost $c_i$ for every time unit spent waiting for service \citep{dai2020queueing,cmu2,cmu3, cmu1}. The problem of minimizing the holding costs is equivalent to the maximization of rewards in a reinforcement learning framework, and we refer to these problems interchangeably in this paper. 
In the queueing networks literature, value-based approximate dynamic programming algorithms have been a popular approach \citep{fluidnetwork, fluid_diffusion, costbound}. 
On the other hand, we aim to use a model-free policy gradient approach to learn the optimal policy directly. Policy gradient methods are favorable in many RL settings because they are effective and scalable to high dimension problems. 

For a single server parallel queueing system, the so-called  $c$-$\mu$ rule is known to minimize the expected holding cost by prioritizing serving the queues based on their holding cost $c_i$ and service time $\mu_i$ \citep{cmu2,cmu3} ratios. Although this rule is simple and easy to implement, it is not obvious  to learn  when the model parameters are unknown. In prior work \citep{cmu1} a greedy algorithm based on the max weight optimization problem was proposed to learn the $c$-$\mu$ for the single server parallel queueing system. 
In this paper, in contrast, we aim to learn the optimal solution directly from the RL environment.
Although the stochastic networks are often far more complex than our simple setting, we believe that our motivations and methods are applicable to more complicated queueing systems, and to other RL applications in general. 

\section{Model description}\label{sec:model}

\subsection{Reinforcement Learning Model}

Reinforcement learning problems are modeled as Markov Decision Processes $\mathcal{M} = (\S, \A, P, r, \mu)$, where $\S$ is  the state space and $\A$ is the action space \citep{bookrl,Sutton2018}. Initial state $s_0$  is chosen  following some distribution $\mu$. At each time step $t = 0, 1, 2, \dots$, the agent takes an action $a_t \in \A$, obtains the immediate reward $r_t = r(s_t, a_t)$, and observes the next state $s_{t+1} $. $P : \S \times \A \to \Delta(\S)$ is a transition function where $P(s'|\, s,a)$ denotes the probability of moving to entering state $s'$ from $s$ after taking action   $a$. 
The state-action record up to time $t$: 
\begin{equation}\tau = (s_0, a_0, r_0, s_1, a_1, r_1,...,s_{t-1}, a_{t-1}, r_{t-1}, s_t)
\end{equation}
is called a trajectory or history.
The average cost of a sample path $\tau$ with horizon $T-1$ is $C(\tau) = \frac{1}{T}\sum_{t=0}^{T-1}  c(s_{t}, a_{t}).$
Our goal is to find a differentiable parameterized policy function $\pi(A)$ that minimizes the expected long term average cost:
\begin{equation}\label{eq:std_prob}
\min_{A} \left\{ J(A) := \E_{\tau \sim \pi(A)}{C(\tau) = \E_{\tau \sim \pi(A)} \left[\lim_{T \to \infty}\frac{1}{T}\sum_{t=0}^{T-1}  c(s_{t}, a_{t})\right]} \right\},
\end{equation}
where $C(\tau)$ is the respective average cost for trajectory $\tau$.

\subsection{Formulation of Parallel Queue System as a Reinforcement Learning Environment}
\begin{figure}[H] 
\centering
\includegraphics[height=0.135\textheight]{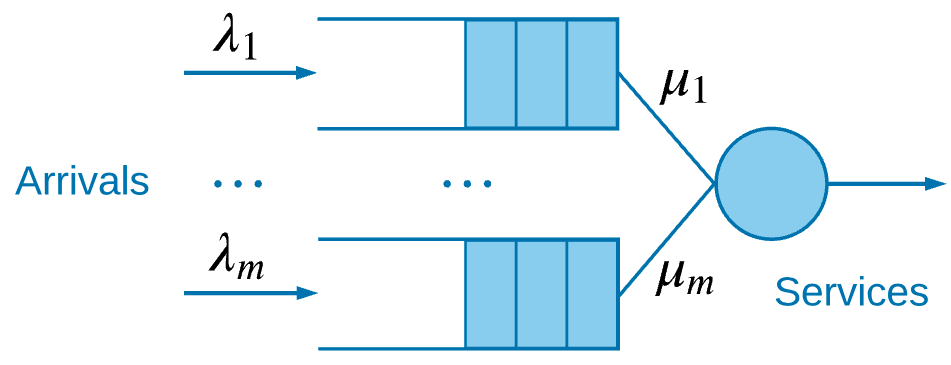}
\caption{Simple parallel queueing system with one service unit. }
\label{fig_sys}
\end{figure}

We consider a parallel system with $m$ queues. This model is simple and has been investigated thoroughly in classical queueing literature \citep{cmu2},
thus, it is a good reference example to develop a better understanding of RL methods. We  assume that the jobs of class $i$ arrive to the system following Poisson processes with respective rates $\lambda_i$, $i = 1, \dots, m$.
There is one server which processes one job at a time. Let us assume that the processing times for class $i$ jobs are i.i.d., having exponential distribution with the respective service rates $\mu_1, \mu_2, \dots, \mu_m $ (thus the expected processing time for a job in class $i$ is $\frac{1}{\mu_i}$. Our model is described in Figure \ref{fig_sys}.
The corresponding load condition for this system is: 
\begin{align}\label{eq_load}
    \rho=\frac{\lambda_1}{\mu_1} + \frac{\lambda_2}{\mu_2} + \dots + \frac{\lambda_m}{\mu_m} < 1
\end{align}
We assume that a decision time occurs when a new job arrives to the system or when the server completes a service. We also use a preemptive policy - when a decision occurs and the server is still in the middle of processing, then the server preempts the unfinished class and follows the new action guided by the decision. Since the inter-arrival times and service times are exponentially distributed, this preemptive process does not alter this assumption. 
The system state is captured in a vector $S = (s,p) \in \N^{m}\times \N$, where $s\in \N^{m}$ is the observable state that shows the number of jobs from each class $1, \dots, m$ waiting in the system, and the last coordinate $p \in \{1,\dots, m\}$ shows the server location $p$ while it is processing a job of class $p$. While a queueing system may have infinite capacity in classic theoretical settings, it is more reasonable to consider a finite buffer in real-life applications. We will consider the latter case in our motivations and empirical findings, where every queue in the system has a capacity of $N$ jobs, and any new arrival of a class which finds $N$ jobs in that class will be lost (without affecting the cost). 

Somewhat contrary to intuition, we allow the server to idle even when there is a job in the system, which implies that there is a total $m + 1$ 
-  actions $1=1,\ldots, m$ select job for class $i$ for processing and the last action lets the server be idle. The server \textbf{always} chooses the next position based on the action provided by the policy and may choose to to serve the queue (class) with no 
job, or it may choose to remain idle even when there is a job in the system. A preemptive policy and an idle action allow us to model the problem more efficiently, and they have been common practices in queueing literature \citep{preemptive,preemtive2,dai2020queueing}. Since the state transition times are exponentially distributed, we can apply the uniformization technique to this problem \citep{unif_cont_discrete, mdp_uniformization}. For the queueing problem after uniformization, the new decision times are determined by the arrival times of a Poisson process with uniform rate that is independent of the current underlying state. 
Classical queueing theory characterizes the existence of a stationary Markovian policy when the load condition \eqref{eq_load} is satisfied (e.g. \citep{stability}).

We let $c \in \R^m$ be the holding cost vector, and let the holding cost be $c^\top s$ where $s \in \R^m$ is the observable state. In the RL setting, the corresponding reward function is $-c^\top s$. The following Theorem characterizes the optimal policy for this system.  

\begin{thm}
For a parallel single server queueing system with infinite buffer, the optimal policy is the {\em priority policy} based on the $c$-$\mu$ rule: The server selects the job in queue $i^*=\argmax \{ c_i \mu_i|\, \text{for\ all\ } i \text{\ such\ that\ } s_i > 0\}$ (i.e. choose the job in the class with the largest cost - expected processing time  ratio among classes that have jobs waiting).
\end{thm}

This theorem is a classical result \citep{cmu2,cmu3, cmu1}.
We denote the priority policy by $\pi_P$. It is important to note that we can only guarantee the optimality of $\pi_P$ under the condition that the system has infinite length and the load condition is satisfied. When the capacity $N$ is big enough, it is natural to expect that the system dynamic is not very different from the infinite buffer case and the priority policy still approximates the optimal 
solution. In the next section, we discuss the parameterization choices and how to represent the priority policy $\pi_P$ in these settings.

\section{New parameterization for queueing models}\label{sec:alg}

\subsection{Parameterizations for Policy Optimization}

We consider a simple class of linear parameterized policies  that map from the state space to a distribution over the action space. Given the state vector $s\in \N^{m}$ a policy is defined matrix a $A \in \R^{(m+1)\times m} $ as follows
\begin{align}
    \pi(A) =  \text{softmax} (z) = \text{softmax} (As),
\end{align}
where $z = As \in \R^{m+1}$ is the linear function of the observable state vector and 
the \textit{softmax} function  defined as 
\begin{align}
\text{softmax}(z) = \left(\frac{e^{z_1}}{ \sum_{k=1}^{m+1} e^{z_k}}; \dots ; \frac{e^{z_{m+1}}}{ \sum_{k=1}^{m+1} e^{z_k}} \right)^\top \in \R^{m+1}. \end{align}
Under such a policy, the probability of any action for any state vector is a number between $0$ and $1$. In contrast, the priority policy $\pi_P$ is a threshold policy, that is given any state vector the policy chooses one action with probability one. This policy cannot be represented as a differentiable function and thus is not amenable to gradient based algorithms.  Our class policy is differentiable but  does not contain $\pi_P$. In the next section we show that $\pi_P$ can be approximated arbitrarily closely with our class of parameterized policies.

\subsection{Representation of Priority Queue}\label{subsec:optimality}
Without loss of generality,  let us assume that  $c_1\mu_1 \geq c_2 \mu_2\geq  \dots \geq c_m \mu_m $. 
In other words, the deterministic priority policy $\pi_P$ chooses to serve at queue $i$ if and only if $i = \arg\min \{i: s_i > 0\}$ where $s\in \N^m$ is the current observable state vector. The following Theorem shows a trajectory that follows the priority policy. 
\begin{thm}\label{thm_1}
Consider the call of  linear policies  described in Section 3.1 and assume that $s_i \leq Q$ for every $i =1, \dots, m$ for a given state vector. Let $\{A_k, k \geq 0\}$ be the sequence 
\begin{align}
     A_k =  \begin{bmatrix}
    (Q+1)^m \cdot k + 1&\dots &1\\
    \dots &(Q+1)^i \cdot k + 1&\dots\\
    1 &\dots & (Q+1) \cdot k + 1\\
    1 & \dots & 1
    \end{bmatrix} \in \R^{(m+1)\times m},
\end{align}
and $\pi_k = \pi(A_k) \in \R^{m+1}$ be the policy corresponding to the matrix $A_k$. We then have
\begin{itemize}
\item The starting point of the sequence $\pi_0$ is a random policy that chooses every action with equal probability. 
\item The sequence $\pi_k$ converges to the priority policy: $\pi_k \to \pi_{P}$ when $k \to \infty$. 
\item There does not exist a bounded sequence of matrices $A_k^\prime$ such that  sequence $ \pi(A_k^\prime) $ converges to the priority policy: $\pi_k \to \pi_{P}$ when $k \to \infty$. 

\end{itemize}

\end{thm}

This theorem suggests some interesting observations about the priority policy $\pi_P$ and our linear policy class. It shows the existence of a  trajectory which starts at a random policy and eventually approximates  priority policy $\pi_P$ arbitrarily closely. This implies that the linear class is sufficient  to  approximately learn the priority policy. Though there may be different trajectories in the parameter space that leads to the policy $\pi_P$, Theorem \ref{thm_1} shows  that the in order to converge to the priority sequence, the policy parameters have to grow arbitrarily large, moreover, their growth may be as fast as exponential in the number of queues. 
We note that this result, with some possible variations, is applicable in many other settings, where  optimal control policies are deterministic and can be similarly expressed as a limit of stochastic policy sequences parameterized by the softmax function.

Based on these observations, we consider an alternative parameterization for $\pi$. We let $B = \ln (A)$ where $\ln$ is an element-wise operator and $A$ is a matrix with positive elements. Let $B$ be the parameters of policy $\bar{\pi}$, we let the  parameterization be
\begin{align}
    \bar{\pi}(B) =  \text{softmax} (e^B s),
\end{align}
where the softmax function is defined above. We denote this scheme as the logarithm-scale parameterization of policy $\bar{\pi}$. It follows immediately that $\pi(A) = \bar{\pi}(B)$ for every matrix $A, B$ satisfying that $A = e^B$.
The sequence of matrices $B_k$ associated with $A_k$ defined in  Theorem \ref{thm_1} is: 
\begin{align}
    B_k &=  \begin{bmatrix}
    \ln[(Q+1)^m \cdot k + 1]&\dots &0\\
    \dots &\ln[(Q+1)^i \cdot k + 1]&\dots\\
    0 &\dots & \ln[(Q+1) \cdot k + 1]\\
    0 & \dots & 0
    \end{bmatrix} \in \R^{(m+1)\times m}.\nonumber
\end{align}
The diagonal elements $B_k$ also grow infinitely large but much slower than those of  $A_k$,  which may result in  better behavior  of the optimization problem (e.g. smoother  objective) when using this parameterization. The question arises if the logarithmically parametrized policy class is more restrictive than the linear class, because it generates only positive-element matrices. The following Theorem suggests that there is an equivalent positive matrix policy to every matrix $A$.

\begin{thm}\label{thm_2}
Let us consider the standard linear parameterization. 
For every matrix A, there exists a  matrix $A'$ that $\pi (A') = \pi (A)$ and $A' > 0$.
\end{thm}

\section{Gradient Estimators for Policy Optimization}\label{sec:estimators}
 We now turn to describing the key ingredient in policy optimization - the function and gradient estimates. 
The function value estimator, which we call {\em inexact zeroth-order oracle} is average cost function induced by the parameterized policy $\pi(A)$. Let $n$ be the number of sample paths and $T-1$ be the time horizon, we define the zeroth-order oracle $\tilde{J}(A)$ as
\begin{align}
    \tilde{J}(A)  = \frac{1}{n} \sum_{i = 1}^n C(\tau_i) = \frac{1}{n} \sum_{i = 1}^n  \left[\frac{1}{T}\sum_{t=0}^{T-1}  c(s_{i,t}, a_{i,t})\right],
\end{align}
where $\tau_i$ is the $i$-th sample path following policy $\pi(A)$. 

We now describe two popular ways of estimating gradients, thus providing {\em inexact first-order} oracles. 

\textbf{Finite Difference Estimator.}
The Finite Difference (FD) estimator $\widetilde{\nabla} J_1 (A)$ has the form 
\begin{align}
      \widetilde{\nabla} J_1(A)  = \frac{1}{M} \sum_{i=1}^M \frac{\tilde{J}(A + u \cdot \vec{v_i}) - \tilde{J}(A) }{u} \vec{v_i},
\end{align} 
where $\vec{v_i}, i=1, \dots, M$ is a set of vectors and $u$ is the finite difference step parameter. 
The essential of this gradient estimator comes from the difference between the (noisy) function values at a parameter $A$ and the reference point $A + u \cdot \vec{v_i} $. 
In this paper, we choose the canonical setting where $\vec{v_i}, i = 1, \dots, M$ are the unit vectors of the parameter space \citep{fd_pde,aaai_fd2}. 
On the other hand, the standard Gaussian random vector is another popular choice for the increment direction \citep{nips_fd1}.


\textbf{Policy Gradient Estimator.}
The classical REINFORCE (PG) method uses the likelihood ratio function to derive its gradient estimator \citep{stoc_opt, glynn,Williams1992}. This estimator is based on the policy gradient theorem \citep{Sutton1999} showing that
$\nabla J(A) = \E_{\tau \sim \pi(A)}\left[ \nabla \log p_A(\tau) R(\tau)\right]$,
where $p_A$ denotes the probability of trajectory $\tau$ under policy $\pi(A)$.
However, the likelihood ratio estimator usually has high variance and REINFORCE estimator also suffers from that weakness. 
Variance reductions techniques have shown improvements to improve the efficiency of this estimator \citep{Sutton2018, Zhao2011, Baxter2001,Kakade2001}. 
In this paper, we choose to present the common practice that adds a baseline term to the final estimator \citep{Zhao2011, pg_robotics,pg_cartpole}. For the time horizon $T$ and  number of sample paths $N$, our Policy Gradient (PG) estimator  $\widetilde{\nabla} J_2 (A)$ is computed as
\begin{align}
      \widetilde{\nabla} J_2(A)  &= \frac{1}{N} \sum_{i=1}^N \left[\frac{1}{T}\sum_{t=0}^{T-1} \nabla_A  \log \pi(A)(s_{i,t}, a_{i,t})\right] \left[\frac{1}{T}\sum_{t=0}^{T-1}  r(s_{i,t}, a_{i,t}) - b \right], \text{ where }\nonumber\\
      b   &= \frac
          {\frac{1}{N} \sum_{i=1}^N \left[\frac{1}{T}\sum_{t=0}^{T-1} \nabla_A  \log \pi(A)(s_{i,t}, a_{i,t})\right]^2 \left[\frac{1}{T}\sum_{t=0}^{T-1}  r(s_{i,t}, a_{i,t})\right]}
          {\frac{1}{N} \sum_{i=1}^N \left[\frac{1}{T}\sum_{t=0}^{T-1} \nabla_A  \log \pi(A)(s_{i,t}, a_{i,t})\right]^2}.
\end{align} 

\begin{figure}[H] 
\centering
\includegraphics[height=0.165\textheight]{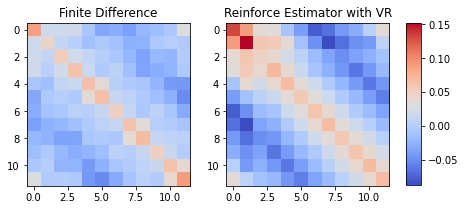}
\includegraphics[height=0.165\textheight]{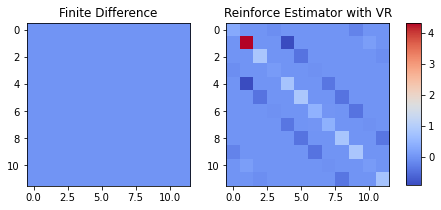}
\caption{Comparisons of the covariance matrices of Finite Difference (FD) estimator and REINFORCE (PG) estimator with Variance Reduction, computed at random policy (left) and (optimal) priority policy (right).  For the random policy, PG estimator behaves in the relatively same scale as FD estimator. However, FD is still significantly better than PG estimator with variance reduction at the priority policy.}
\label{fig_cov}
\end{figure}

\textbf{Comparing two gradient estimators}. In order to further understand the properties of our first-order oracles, we conduct a small experiment to investigate their behavior. Both of these gradient estimators have been popular in the literature, however,  as far as we know they have not been compared side by side. The main reason for this is that they are typically implemented within algorithms customized the the specific choice of first order oracle.  
In our experiment, we choose the number of samples $n$ and $N$ such that two estimators have the same sample complexity. We compute two estimators at two reference points: at a random policy and at the (optimal) priority policy. We evaluate the covariance matrices of these stochastic estimators and plot the results in Figure \ref{fig_cov}.

This experiment shows that the REINFORCE estimator has similar or slightly higher variance than Finite Difference  estimator at the random policy. When we compute the two estimators at the (optimal) priority policy, the FD estimator has an extremely low variance and outperforms the PG estimator. Thus, we conclude that the FD estimator is more reliable than PG estimator when the algorithm reaches the neighborhood of minimizers.
\section{Numerical Experiments}\label{sec_experiment}

In this section, we experiment with two parameterizations: the standard linear parameterization and the logarithm-scale parameterization proposed in Section \ref{subsec:optimality}. We implement two gradient estimators and employ an adaptive line-search algorithm (ALOE) to optimize the cost function. Unlike other more complex policy optimization algorithms (e.g. TRPO and PPO \citep{Schulman2015,Schulman2017}) this method, describe it in Algorithm \ref{alg_aloe}, has theoretical complexity guarantees for both types on first order oracles we describe above and allows for bias in zeroth-order and first-order oracles  \citep{billy_miaolan_katya, katya_linesearch}.

\begin{algorithm}[hpt!]
  \caption{ Adaptive Line-search with Oracle Estimations (ALOE) \citep{ katya_linesearch}
}
  \label{alg_aloe}
\begin{algorithmic}
  \STATE {\bfseries Input:} Parameter $\epsilon_f$ of the zeroth-order oracle, starting point $A_0$, max step size $\alpha_{\max} > 0$, initial step size $\alpha_{0} < \alpha_{\max}$, constants $\theta, \gamma \in (0, 1)$. 
  \FOR{$k=0,1,2, \cdots$}
  \STATE Compute gradient approximation $ g_k = \widetilde{\nabla} J(A_k)$ 
  \STATE Check sufficient decrease: 
  \STATE \hspace{1.5em} Let $A_k^+ = A_k - \alpha_k g_k$. Generate the zeroth-order oracle $\tilde{J}(A_k)$ and $\tilde{J}(A_k^+)$. 
  \STATE \hspace{1.5em} Check the modified Armijo condition:
    \begin{align}\label{eq_conditions}
        \tilde{J}(A_k^+) \leq \tilde{J}(A_k) - \alpha_k \theta \|g_k\|^2 + 2 \epsilon_f.
    \end{align}
  \STATE \textbf{\textit{Successful step:}} If \eqref{eq_conditions} holds, then set $A_{k+1} = A_k^+$ and $\alpha_{k+1} = \min \{ \alpha_{\max}, \gamma^{-1} \alpha_k \}$. 
  \STATE \textbf{\textit{Unsuccessful step:}} Otherwise, set $A_{k+1} = A_k$ and $\alpha_{k+1} =  \gamma \alpha_k$. 
  \ENDFOR
\end{algorithmic}
\end{algorithm} 

\begin{figure}[H] 
\centering

\includegraphics[height=0.15\textheight]{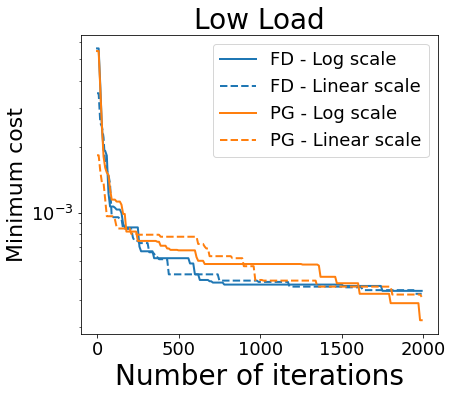}
\includegraphics[height=0.15\textheight]{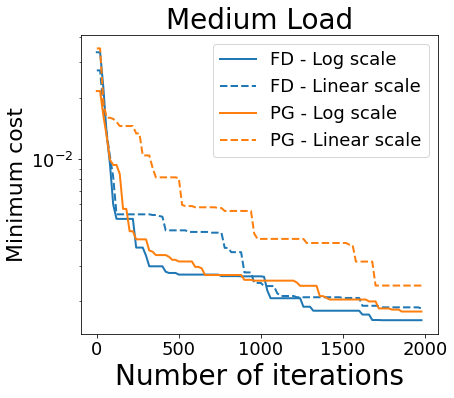}
\includegraphics[height=0.15\textheight]{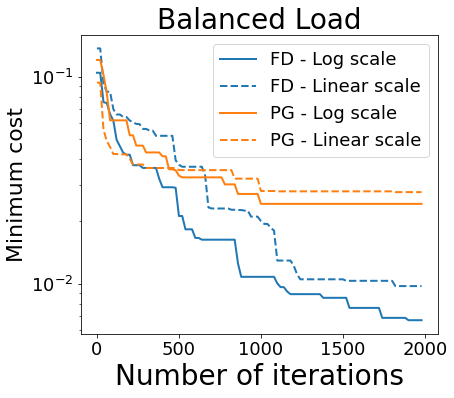}
\includegraphics[height=0.15\textheight]{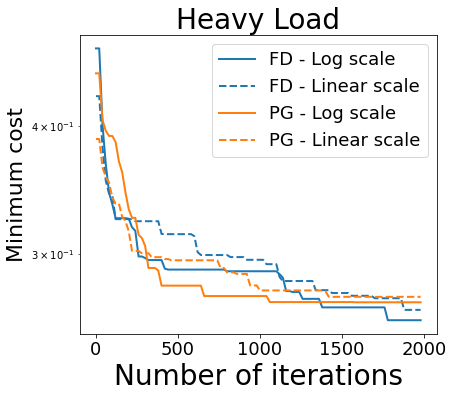}
\caption{Comparisons of the minimum cost function achieved by iteration, for two parameterization scales: logarithm (proposed scale, presented in continuous solid lines) and linear scale (presented in dashed line). We show the performance of Finite Difference (FD) estimator and REINFORCE estimator with variance reduction (PG) in different colors. 
 }
\label{fig_cost}
\end{figure}

\textbf{Compare the cost function for two parameterizations.}
We experiment with two parameterizations of the policy with various load conditions of the system from light to heavy traffic. We summarize our setting in Table \ref{tab_reward} and highlight the best results obtained by our algorithms there. In addition, we plot the training progress in Figure \ref{fig_cost}. 
In each setting, we tune the algorithms using grid search and pick the best choice of hyper-parameters (e.g. the step sizes $\alpha_{\max}$ and $\alpha_0$) to the training stage. We repeat all the stochastic experiments for 5 random seeds, then report the average performance. 
We present the detailed implementation settings in the Appendix. 

Table \ref{tab_reward} shows that in most cases, the logarithm-scales perform better than the standard linear-scales. The Finite Difference estimator often shows better performance than Policy Gradient estimator, from medium to heavy load conditions.
We also report the confidence interval (CI) of the cost function evaluated at the optimal (priority) policy. 
Since ALOE accepts steps based on a modified function reduction condition, the obtained cost values tend to be closer to the lower end of the corresponding confidence intervals. 

\renewcommand{\arraystretch}{1.25}
\begin{table}[H]
    \centering
    \caption{Summary of the load conditions in our experiment setting, and the best performance of each  parameterization (Linear scale and  Logarithm scale). The load parameter $\rho$ is presented in equation \eqref{eq_load}. The cost function is normalized to $[0,1]$.
    We report the best cost function achieved in each setting and note the first order oracle that yields the best result. The last row reports the confidence intervals of the cost function at the optimal (priority) policy in each setting.\\}
    \label{tab_reward}
\begin{tabular}{|l|l|l|l|l|}
    \hline
    Settings & Low Load & Medium Load & Balanced Load & Heavy Load\\
    \hline
    Service rates & $(18,9,6)$ & $(9,4.5,3)$ & $(6,3,2)$ & $(3,2,1)$\\
    \hline
    Load parameter $\rho$& $1/3$ &  $2/3$ & $1$ &  $11/6$\\
    \hline
    Linear scale & 0.0003 (PG) &0.0018 (FD) & 0.0097 (FD) & 0.2647 (FD)\\
     \hline
    Logarithm scale & \textcolor{blue}{0.0004 (PG) }& \textcolor{blue}{0.0016 (FD) }& \textcolor{blue}{0.0066 (FD)} &\textcolor{blue}{0.2586 (FD)}\\
     \hline
     Optimal cost CI&$[0.0002, 0.0015]$&$[0.0012 ,0.0032]$&$[0.0079, 0.0093]$&$[0.2611, 0.4037]$\\
     \hline
\end{tabular}
\end{table}

\renewcommand{\arraystretch}{1.2}

\begin{figure}[H] 
\centering
\includegraphics[height=0.15\textheight]{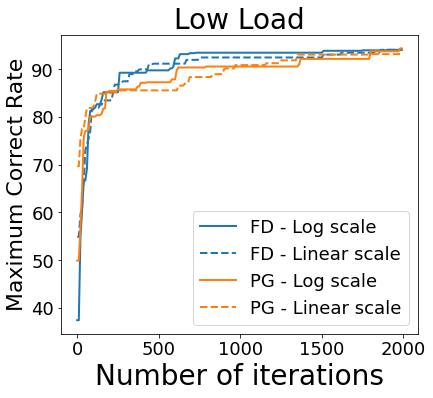}
\includegraphics[height=0.15\textheight]{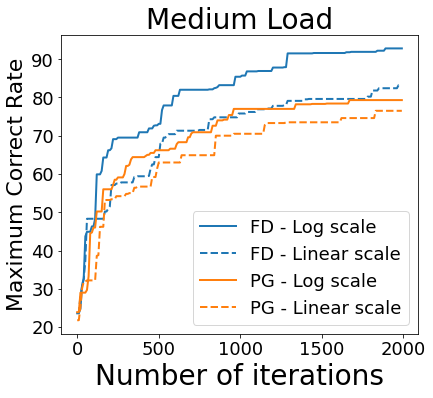}
\includegraphics[height=0.15\textheight]{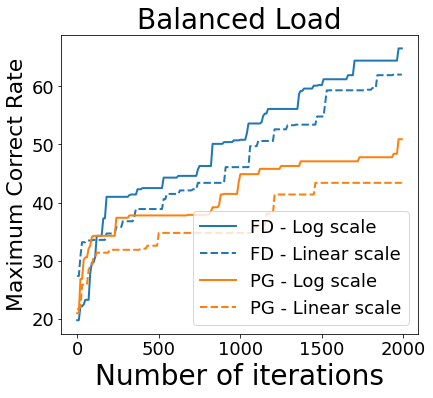}
\includegraphics[height=0.15\textheight]{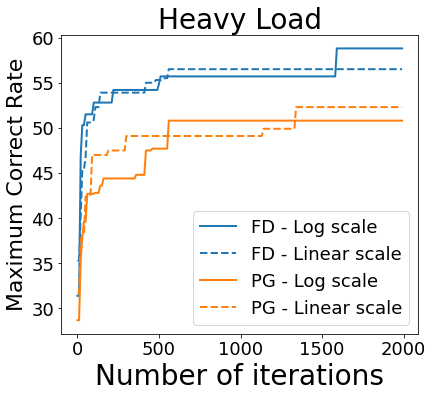}
\caption{Comparisons of the maximum correct rate (in percent) achieved by iteration for two parameterization scales: logarithm (proposed scale, presented in continuous solid lines) and linear scale (presented in dashed line). We show the performance of Finite Difference (FD) estimator and REINFORCE estimator with variance reduction (PG) in different colors. 
}
\label{fig_correct_rate}
\end{figure}

\textbf{Compare the learned policy with the priority policy.}
Now we compare the learned policies with the priority policy and show how they behave relatively to the desired outcome. For every pair of state and action induced by policy $\pi(A)$, we compute the rate that $\pi(A)$ chooses the "correct" action dictated by priority policy $\pi_P$. We keep measuring that correct rate along the training process, and report the results in Figure \ref{fig_correct_rate}. 
Our summary in Table \ref{tab_correct} shows that the logarithm scale usually has a better performance than the normal linear parameterization. In addition, it suggests that the algorithms learn the correct action more easily in the low and medium load conditions, with higher than 90 percent accuracy rates. In the balanced and heavy load settings, it is more difficult for the algorithms to learn the correct actions, however, it still achieves an accuracy rate of approximately 60 percent. We present the experiment setting in the Appendix.

\begin{table}[H]
    \centering
    \caption{The best correct rate in each setting achieved by our algorithm, with the corresponding gradient estimators that yields the best reported result.\\}
    \label{tab_correct}
    \begin{tabular}
    {|l|m{7em}|m{7em}|m{7em}|m{7em}|} 
    \hline
    Settings & Low Load & Medium Load & Balanced Load & Heavy Load\\
    \hline
    Service rates & $(18,9,6)$ & $(9,4.5,3)$ & $(6,3,2)$ & $(3,2,1)$\\
    \hline
    Load parameter $\rho$ & $1/3$ &  $2/3$ & $1$ &  $11/6$\\
    \hline
    Linear scale  & 94.3 (FD) & 83.2 (FD) & 62.0 (FD) & 56.5 (FD)\\
     \hline
    Logarithm scale & \textcolor{blue}{94.4 (PG)} & \textcolor{blue}{92.8 (FD)}& \textcolor{blue}{66.5 (FD)} & \textcolor{blue}{58.8 (FD)}\\
     \hline
    \end{tabular}
\end{table}
\section{Conclusion}\label{sec_conclusion}
We investigate the optimization aspect of queueing reinforcement learning environments, by proposing new logarithm-scale parameterization  and also comparing two gradient estimators side-by side within a stochastic step search algorithm. In most settings, the logarithm-scale parameterization and the finite difference based gradient estimator  show better performance than the linear parametrization and policy gradient estimators. 

\section*{Acknowledgements}
The authors sincerely thank 
Jamol J. Pender and Shane G. Henderson for their insightful discussions and useful suggestions to complete this project. The work of Trang H. Tran and Katya Scheinberg have partly been supported by the ONR Grant N00014-22-1-2154.


\appendix

\appendix
\newpage

\section*{APPENDIX}

\section{Technical Proofs}
\subsection*{Proof of Theorem \ref{thm_1}}

\begin{proof}
From the proof of Theorem \ref{thm_2}, we note that it is sufficient to prove the statements for the following sequence: 
\begin{align*}
     A''_k =  \begin{bmatrix}
    (Q+1)^m \cdot k &\dots &0\\
    \dots &(Q+1)^i \cdot k&\dots\\
    0 &\dots & (Q+1) \cdot k\\
    0 & \dots & 0
    \end{bmatrix} \in \R^{(m+1)\times m},
\end{align*}
because $A_k = A''_k +1$ and $\pi(A_k) = \pi(A''_k)$. 

The first statement of Theorem \ref{thm_1} when $k=0$ is straightforward. We note that $A''_0 = 0$ and the softmax function at zero vector returns a random policy that chooses every action with equal probability. 

Now we move to the second statement. Let $z_k = A''_k s$, we have: 
\begin{align*}
     z_k &=  \left(
    (Q+1)^m \cdot k s_1; \dots; (Q+1)^i \cdot k s_{m-i} ;\dots; (Q+1) \cdot k s_m ;
    0\right)\\
    &= k \left(
    (Q+1)^m  s_1; \dots; (Q+1)^i  s_{m-i} ;\dots; (Q+1)  s_m ;
    0\right) = k\cdot y,
\end{align*}
where $s_i$ is the $i$-th element of the observable state $s$, and 
\begin{align*}
    y &= \left(
    (Q+1)^m  s_1; \dots; (Q+1)^i  s_{m-i} ;\dots; (Q+1)  s_m ;
    0\right)
    \in \R^{m+1},
\end{align*}

Now let  $j = \arg\min \{j: s_j > 0\}$ and assuming that $j$ is not an idle action, we have that the deterministic priority policy $\pi_P$ chooses to serve at queue $j$. We prove that $y_j$ is the largest element of vector $y \in \R^{m+1}$. 

Firstly, we have $s_i = 0$ for every $i < j$ and therefore $y_i = 0$  for every $i < j$ and $y_j = (Q+1)^{m-j}  s_{j} > 0 = y_i $ for every $i < j$. In addition, since $s_i \leq Q$ for every $i > j$, we have that $s_i < Q+1$ and 
$$y_i = (Q+1)^{m-i}  s_i < (Q+1)^{m-i+1} \leq (Q+1)^{m-j} \leq (Q+1)^{m-j}  s_{j} = y_j.$$

Hence $y_j$ is the largest element of vector $y \in \R^{m+1}$. Now we consider the softmax policy induced by $z_k = k\cdot y$: 
$$\text{softmax}(k\cdot y) = \left(\frac{e^{k y_1}}{ \sum_{i=1}^{m+1} e^{k y_i}}; \dots ; \frac{e^{k y_{m+1}}}{ \sum_{i=1}^{m+1} e^{k y_i}} \right)^\top \in \R^{m+1}. $$

Note that $k y_j$ is the largest element of vector $k y \in \R^{m+1}$. Hence when $k$ goes to infinity, the vector $\text{softmax}(k\cdot y)$ converges to the $j$-th unit vector in $\R^{m+1}$. This is the deterministic priority policy $\pi_P$ chooses to serve at queue $j$.

Now we move to the final statement. Let us fixed the state $s$ we have that $(\pi_P)_i = 0$ for every incorrect action, i.e. action $i\neq j$. We assume the contradiction that there exist a sequence $A'_k$  such that all the (absolute value) of elements of $A'_k$ are upper bounded by a constant $R$, and  sequence $ \pi(A_k^\prime) $ converges to the priority policy: $\pi_k \to \pi_{P}$ when $k \to \infty$.

It is easily seen that all the (absolute value) of elements of the vector $z'_k = A'_k s$ are upper bounded by $m R Q$, since the state space is bounded by $Q$. Thus the denominator of the softmax vector (i.e. sum of vector $\exp({z'_k})$) is upper bounded by $(m+1) e^{mRQ}$. However, for an action $i \neq j$, we have the softmax value converge to 0. This shows that $\exp({(z'_k)_i})$ converges to 0 and $(z'_k)_i \to -\infty$, which contradicts to the fact that $z'_k$ and $ A'_k$ are upper bounded.

Hence we complete the proof. Note that the construction of sequence $A_k$ is not unique, and there may be different trajectories in the parameter space that leads to the policy $\pi_P$. However, this result suggests that the limit of such trajectory need to have the ability to distinguish (and let dominate) the more important queue. Such ability leads to the big elements of $A_k$ when $k$ is large. 
\end{proof}
\newpage
\subsection*{Proof of Theorem \ref{thm_2}}

Let us consider the standard linear parameterization. 
For every matrix $A$, there exists a shifted matrix $A'$ that $\pi (A') = \pi (A)$ and $A' > 0$.
\begin{proof}
We recall that each policy is represented by a matrix $A \in \R^{(m+1)\times m} $, and the output distribution is 
\begin{align*}
    \pi(A) =  \text{softmax} (z) = \text{softmax} (As),
\end{align*}
where $z = As \in \R^{m+1}$ is the linear function output.
The \textit{softmax} function is defined as $$\text{softmax}(z) = \left(\frac{e^{z_1}}{ \sum_{k=1}^{m+1} e^{z_k}}; \dots ; \frac{e^{z_{m+1}}}{ \sum_{k=1}^{m+1} e^{z_k}} \right)^\top \in \R^{m+1}. $$

We prove that for every number $\mu \in \R$, the following statement holds: $\pi (A) = \pi(A + \mu),$ where $A+\mu$ yields the element-wise addition between the elements of $A$ and $\mu$. In fact, we have 
\begin{align*}
    \pi(A+\mu) &= \text{softmax} ((A+\mu)s) \\
    &= \text{softmax} \left(As + \vec{e}_{m+1}  \cdot \mu \sum_{k=1}^{m+1} s_k\right)\\
    &= \text{softmax} \left(As\right) = \pi(A),
\end{align*}
where $\vec{e}_{m+1}$ is the vector containing all elements 1 in $\R^{m+1}$. The last line follows since the softmax function is invariant with respect to the addition of the scalar of $\vec{e}_{m+1}$.

Finally, choosing $\mu$ such that all the elements of $A + \mu$ is postive, we have $\mu = 1 - \min{A_{ij}}$ and $A' = A + \mu$ satisfying that $\pi (A') = \pi (A)$ and $A' > 0$.
\end{proof}
\section{Experiment Settings}
We describe the detailed implementation settings for our experiments here. 

We model the queueing system dynamics using OpenAI gym package \cite{openaigym}. We implement the zeroth- and first-order oracles for our environment using a warm up stage of 1000 iterations to reduce the initial variance. In addition, the cost function is normalized to $[0,1]$. In more detail, the cost function is $1/(mQ)$ where $m$ is the number of queues in the system and $Q$ is the maximum number of jobs in every queue. Thus, the cost function is 0 when the system has no job, and it is 1 when the system is completely full with $mQ$ jobs. 

The service rates for each setting is given in Table \ref{tab_reward}. We choose $m=3$ and $Q=100$. We choose the number of sample paths $n$ and $N$ such that two estimators have the same sample complexity. The number of sample paths $n$ in the zero-th order oracle is 20. The number of (basis) vector used in the FD estimator is the size of $A$, which is $m(m+1)$. Hence we choose the number of sample paths of 
PG estimator to be $N=[m(m+1) + 1]n$ to match the number of sample paths of 
FD estimator. 

The time horizon $T$ is 100 and we optimize over 2000 iterations of ALOE algorithm. We start all the algorithms at random policy $A=0$. We tune different hyper-parameters for ALOE algorithm using grid search and choose the best parameters in to the final stage. The final setting of ALOE algorithm is: $\epsilon_f = 0.01$, $\alpha_{\max} = 2, \alpha_0 = 0.01$, $\theta = 0.01$ and $\gamma = 0.5$. In Table \ref{tab_reward}, we compute the confidence intervals using 1000 sample paths and with horizon 100. We also use a warm up stage with 1000 iterations.

\newpage

\bibliography{reference}
\bibliographystyle{plainnat}





\end{document}